\documentclass[12pt,a4paper]{article}
\usepackage{latexsym,amssymb,amsmath,mathrsfs}

\usepackage{sectsty}
\usepackage{color}
\usepackage{bm}
\usepackage[T1]{fontenc}
\usepackage[anchorcolor=blue,%
 bookmarks=true,%
 bookmarksnumbered=true,%
]{hyperref}
\usepackage{cite}
\usepackage{pst-all}
\usepackage{graphicx}
\allsectionsfont{\normalsize\sc\center}

\newtheorem{thm}{Theorem}[section]
\newtheorem{prop}[thm]{Proposition}
\newtheorem{cor}[thm]{Corollary}
\newtheorem{lem}[thm]{Lemma}

\newtheorem{remark}[thm]{Remark}

\newenvironment{pf}{\par\begin{trivlist}%
\item[]{\bf Proof.}\ }{\hfill $\square$ \end{trivlist}\par}

\makeatletter \@addtoreset{equation}{section} \makeatother

\newcommand{\R}{\mathbb{R}}

\DeclareMathOperator{\Ric}{Ric}

\renewcommand{\d}{\mathrm{d}}

\newcommand{\s}{\mathfrak{s}}

\newcommand{\eps}{{\varepsilon}}

\title{\large\bf  Rigidity phenomena on lower \boldmath$N$-weighted Ricci curvature bounds with 
\boldmath$\varepsilon$-range for non-symmetric Laplacian}
\author{Kazuhiro Kuwae\thanks{Department of Applied Mathematics, Fukuoka University,
Fukuoka 814-0180, Japan ({\sf kuwae@}  {\sf fukuoka-u.ac.jp}) . Supported in part by JSPS Grant-in-Aid for Scientific Research (KAKENHI) 17H02846 and by fund (No.:185001) from the Central Research Institute of Fukuoka University.}
\ \ and\ \ 
Yohei Sakurai\thanks{Department of Mathematics, Saitama University,
255 Shimo-Okubo, Sakura-ku, Saitama-City, Saitama, 
338-8570, Japan ({\sf ysakurai@rimath.saitama-u.ac.jp}). Supported in part by JSPS Grant-in-Aid for Scientific Research on Innovative Areas ``Discrete Geometric Analysis for Materials Design" 17H06460.
}
}
\date{}
\begin{document}
\maketitle
% \raggedbottom 
\begin{abstract}
Lu-Minguzzi-Ohta \cite{LMO:CompaFinsler} have introduced the notion of a lower $N$-weighted Ricci curvature bound with $\eps$-range,
and derived several comparison geometric estimates from a Laplacian comparison theorem for weighted Laplacian. 
The aim of this paper is to investigate various rigidity phenomena for the equality case of their comparison geometric results.
We will obtain rigidity results concerning the Laplacian comparison theorem, diameter comparisons, and volume comparisons.
We also generalize their works for non-symmetric Laplacian induced from vector field potential.
\end{abstract}


 
%%%%%%%%%%%%%%%%%%%%%
%%%%%%%%%%%%%%%%%%%%%
%%%%%%%%%%%%%%%%%%%%%
\section{Introduction}

%%%%%%%%%%%%%%%%%%%%%
\subsection{Weighted Ricci curvature and comparison geometry}

Let $(M,g,f)$ denote an $n$-dimensional weighted Riemannian manifold, namely, $(M,g)$ is an $n$-dimensional complete Riemannian
manifold, and $f\in C^{\infty}(M)$.
For $N\in]-\infty,\,+\infty]$, the {\it $N$-weighted Ricci curvature} is defined as follows (\cite{BE1}, \cite{Lich}): 
\begin{align*}
\Ric_f^N:=\Ric_g+{\rm \nabla^2}f-\frac{df\otimes df}{N-n},%\label{eq:weightedRici}
\end{align*}
where when $N=+\infty$,
the last term is interpreted as the limit $0$,
and when $N=n$,
we only consider a constant function $f$, and
set $\Ric_f^n:=\Ric_g$.

It is well-known that lower weighted Ricci curvature bounds lead us various comparison geometric results.
In the traditional case of $N\in [n,+\infty[$,
under a curvature condition
\begin{align}
\Ric_f^N\geq Kg \label{eq:RicciLowerBdd}
\end{align}
for $K\in \mathbb{R}$,
such results have been obtained by \cite{Lo}, \cite{Qi}, \cite{WW},
and so on.
On the other hand,
recently,
comparison geometry has begun to be developed in the complementary case of $N\in ]-\infty,n[$ (see e.g., \cite{KL}, \cite{KS}, \cite{Lim}, \cite{LMO:CompaFinsler}, \cite{Mai1}, \cite{Mai2}, \cite{Mineg}, \cite{Oh<0}, \cite{Sak}, \cite{Wy}, \cite{WyYero}).
Wylie-Yeroshkin \cite{WyYero} have introduced a curvature condition
\begin{align}
\Ric_f^1\geq (n-1)\kappa e^{-\frac{4f}{n-1}}g \label{eq:RicciLowerBddWY}
\end{align}
for $\kappa \in \mathbb{R}$,
and presented an optimal Laplacian comparison theorem, Bonnet-Myers theorem, Bishop-Gromov volume comparison theorem.
After that
the first named author and Li \cite{KL} have extended their condition to
\begin{align}
\Ric_f^N\geq (n-N)\kappa e^{-\frac{4f}{n-N}}g\label{eq:RicciLowerBddKL}
\end{align}
with $N\in ]-\infty,1]$,
and also done their comparison theorems.

Lu-Minguzzi-Ohta \cite{LMO:CompaFinsler} have introduced a new curvature condition that interpolates the conditions \eqref{eq:RicciLowerBdd} with $K=(N-1)\kappa$, \eqref{eq:RicciLowerBddWY} and \eqref{eq:RicciLowerBddKL}.
For $N\in ]-\infty,1]\cup [n,+\infty]$,
they have considered the notion of the \emph{$\varepsilon$-range}:
\begin{align}
\varepsilon=0\ \text{ for }\ N=1,\qquad \varepsilon\in ]-\sqrt{\varepsilon_0},\sqrt{\varepsilon_0}[  \ \text{ for }\  N\ne1,n,\qquad \varepsilon\in\R \ \text{ for } \ N=n,\label{eq:epsilonrange}
\end{align}
where
\begin{align*}
\varepsilon_0:=\frac{N-1}{N-n}.
\end{align*}
Here if $N=+\infty$,
then we interpret $\eps_0$ as the limit $1$.
In this range,
they have proposed a curvature condition
\begin{align}
\Ric_f^N\geq c^{-1}\,\kappa\, e^{-\frac{4(1-\varepsilon)f}{n-1}}g\label{eq:LMOLowerBddd}
\end{align}
for $\kappa\in\R$, where $c=c_{N,\varepsilon}\in]0,1]$ is a positive constant defined by 
\begin{align}
c:=\frac{1}{n-1}\left(1-\varepsilon^2 \frac{N-n}{N-1} \right)\label{eq:constant}
\end{align}
if $N\ne1$,
and $c:=(n-1)^{-1}$ if $N=1$. 
When $N\in [n,+\infty[$ and $\varepsilon=1$ with $c=(N-1)^{-1}$, the curvature condition \eqref{eq:LMOLowerBddd} is reduced to \eqref{eq:RicciLowerBdd} with $K=(N-1)\kappa$.
Also,
when $N=1$ and $\varepsilon =0$ with $c=(n-1)^{-1}$, it covers \eqref{eq:RicciLowerBddWY},
and when $N\in ]-\infty,1]$ and $\varepsilon=\varepsilon_0$ with $c=(n-N)^{-1}$, 
it does \eqref{eq:RicciLowerBddKL}. 
Under the condition \eqref{eq:LMOLowerBddd},
they first proved a Laplacian comparison theorem for the distance function,
and derived a diameter bound of Bonnet-Myers type,
and a volume bound of Bishop-Gromov type under density bounds.
Notice that
they have worked on Finsler setting beyond weighted Riemannian setting.

%%%%%%%%%%%%%%%%%%%%%
\subsection{Setting}

The purpose of this paper is to examine rigidity phenomena on the equality case of comparison theorems under the curvature condition \eqref{eq:LMOLowerBddd}.
We work on a general setting such that
a weighted manifold $(M,g,V)$ has a vector field potential $V$ beyond the gradient case of $V=\nabla f$.
Such a weighted manifold is equipped with canonical weighted Laplacian and weighted Ricci curvature.
The {\it weighted Laplacian} is defined by
\begin{equation*}
\Delta_V:=\Delta- g(V, \nabla\cdot),
\end{equation*}
and the {\it $N$-weighted Ricci curvature} is done as follows:
\begin{align*}
{\Ric}_{V}^N:={\Ric}_g+\frac12\mathcal{L}_Vg-\frac{V^*\otimes V^*}{N-n}.
\end{align*}
Here $\mathcal{L}_Vg$ is the Lie derivative of $g$ with respect to $V$,
and $V^*$ denotes its dual $1$-form.  
If $N=n$,
then we always assume that $V$ vanishes such that ${\Ric}_{V}^n={\Ric}_g$.
In the gradient case of $V=\nabla f$,
it coincides with $\Ric_f^N$.

We now describe our setting.
We always fix a point $p\in M$,
and also $N\in ]-\infty,1]\cup [n,+\infty]$ and $\eps \in \mathbb{R}$ in the range \eqref{eq:epsilonrange}.
We consider an arbitrary positive constant $c_p>0$ such that
we choose
\begin{equation}\label{eq:gradcp}
c_p=e^{   -\frac{2(1-\eps)f(p)}{n-1}}
\end{equation}
in the gradient case of $V=\nabla f$.
We further define two lower semi continuous functions $f_{V,p},\,s_{V,p}:M\to \mathbb{R}$ by
\begin{equation*}
f_{V,p}(x):=\inf_{\gamma}\int^{d(p,x)}_0\,g(V,\dot{\gamma}(\xi))\,\d \xi,\quad s_{V,p}(x):=c_p\inf_{\gamma}\int^{d(p,x)}_0\,e^{      -\frac{2(1-\eps)f_{V,p}(\gamma(\xi))}{n-1}       }\,\d \xi,
\end{equation*}
where the infimum is taken over all unit speed minimal geodesics $\gamma:[0,d(p,x)]\to M$ from $p$ to $x$.
In the gradient case of $V=\nabla f$,
we see $f_{V,p}(x)%:
=f(x)-f(p)$.
Furthermore,
$s_{V,p}$ is called the {\it re-parametrized distance} from $p$ (cf. \cite{WyYero}).
For a continuous function $\kappa:[0,+\infty[\to \R$,
we also define a function $\kappa_{V,p}:M\to \mathbb{R}$ by
\begin{equation*}
\kappa_{V,p}:=\kappa \circ s_{V,p}.
\end{equation*}
Our setting is as follows:
\begin{align}
\Ric_V^N\geq c^{-1}\,c^2_p\,\kappa_{V,p}\, e^{-\frac{4(1-\varepsilon)f_{V,p}}{n-1}}g,\label{eq:WeightedLowerBddd}
\end{align}
where $c$ is defined as \eqref{eq:constant}.
In the gradient case of $V=\nabla f$, and the case where $\kappa$ is constant,
this is reduced to \eqref{eq:LMOLowerBddd}.
Under the condition \eqref{eq:WeightedLowerBddd} with non-gradient potential,
for instance,
Bakry-Qian \cite{BQ2} have obtained a Laplacian comparison theorem and a volume comparison theorem of Bishop-Gromov type for invariant measures when $N\in [n,+\infty[$ and $\varepsilon=1$ with $c=(N-1)^{-1}$,
Kuwada \cite{K} has given a diameter comparison theorem of Bonnet-Myers type and a maximal diameter theorem of Cheng type when $N\in [n,+\infty[$ and $\varepsilon=1$ and $\kappa$ is a positive constant,
Wylie \cite{Wy} has proven a splitting theorem of Cheeger-Gromoll type when $N\in ]-\infty,1]$ and $\kappa =0$,
and the first named author and Shukuri \cite{KS} have studied various comparison geometric properties when $N\in ]-\infty,1]$ and $\varepsilon=\varepsilon_0$ with $c=(n-N)^{-1}$.
We also aim to develop such comparison geometry.

%%%%%%%%%%%%%%%%%%%%%
\subsection{Organization}

In Section \ref{sec:Laplacian},
we will produce a Laplacian comparison theorem,
and its rigidity properties for the equality case,
which is a key ingredient of the proof of our rigidity theorems (see Theorem \ref{lem:Laplacian comparison} and Lemma \ref{lem:LaplacianRigidity}).
Our Laplacian comparison is a generalization of the one that has been obtained by Lu-Minguzzi-Ohta \cite{LMO:CompaFinsler} in the gradient case.
In \cite{LMO:CompaFinsler},
they have derived it from the so-called Bishop inequality,
which is based on an algebraic calculation.
To analyze the rigidity phenomena in more detail,
inspired by the original work of Wylie-Yeroshkin \cite{WyYero},
we will deduce it from the Bochner formula.

In Section \ref{sec:Diameter},
we study diameter comparison theorems of Bonnet-Myers type,
and maximal diameter theorems of Cheng type for the equality cases.
We will obtain two maximal diameter theorems.
The first one is a generalization of the one that has been proved by Wylie-Yeroshkin \cite{WyYero} under the curvature condition \eqref{eq:RicciLowerBddWY} concerning the diameter of a conformally deformed Riemannian metric (see Theorem \ref{thm:diameter rigidity} and Corollary \ref{cor:corodiameter rigidity}).
The second one is a new result even in the setting of Wylie-Yeroshkin \cite{WyYero}.
We will assume not only the curvature bound but also a density bound %the  
in the gradient case.
We characterize the equality case of a diameter comparison by standard sphere with constant density (see Theorem \ref{thm:finite diameter rigidity} and Corollary \ref{cor:cor finite diameter rigidity}).

In Section \ref{sec:Volume},
we investigate volume comparisons.
We first deduce an absolute comparison of Bishop type,
and a relative volume comparison of Bishop-Gromov type concerning a weighted volume of sub-level sets of the re-parametrized distance that have been obtained in the setting of Wylie-Yeroshkin \cite{WyYero} (see Propositions \ref{thm:absolute volume comparison} and \ref{thm:relative volume comparison}).
We also establish a rigidity theorem for the equality case of them (see Theorem \ref{thm:volume growth rigidity}).

In Section \ref{sec:Compactness},
we examine some compactness properties in our setting.
First,
we consider the notion of $\eps$-completeness from which we derive a compactness property (see Proposition \ref{prop:compact}).
We further show a theorem of Ambrose type (see Theorem \ref{thm:AmbroseMyers}).
%%%%%%%%%%%%%%%%%%%%%
%%%%%%%%%%%%%%%%%%%%%
%%%%%%%%%%%%%%%%%%%%%
\section{Laplacian}\label{sec:Laplacian}

%%%%%%%%%%%%%%%%%%%%%
\subsection{Riccati inequality}

Let $d_{p}:M\to \mathbb{R}$ be the distance function from $p$ defined as $d_{p}(x):=d(p,x)$.
We denote by $U_pM$ the unit tangent sphere at $p$.
We define a function $\tau:U_{p}M \to ]0,+\infty]$ as
\begin{equation*}
\tau(v):=\sup \{t>0 \mid d_{p}(\gamma_{v}(t))=t \},
\end{equation*}
where $\gamma_{v}:[0,+\infty[\to M$ is the unit speed geodesic with $\gamma_v(0)=p$ and $\dot{\gamma}_{v}(0)=v$.

We start from the following Riccati inequality that has been already shown by Lu-Minguzzi-Ohta \cite{LMO:CompaFinsler} in the gradient case of $V=\nabla f$ via algebraic calculation (see the proof of \cite[Proposition 3.5]{LMO:CompaFinsler}).
Inspired by the work of Wylie-Yeroshkin \cite{WyYero},
we give its proof by using Bochner formula (cf. \cite[Lemma 4.1]{WyYero}, \cite[Lemma 3.1]{KS}).
\begin{lem}\label{lem:Riccati}
For all $t \in ]0,\tau(v)[$
we have
\begin{align}\label{eq:Riccati}
&\quad  \left(\bigl(e^{\frac{2(1-\eps)f_{V,p}}{n-1}}\,\Delta_{V}d_{p}\bigl)(\gamma_{v}(t))\right)'\\ \notag
&\leq      -e^{\frac{2(1-\eps)f_{V,p}(\gamma_{v}(t))}{n-1}}\,\Ric^{N}_{V}(\dot{\gamma}_{v}(t))-c\,e^{-\frac{2(1-\eps)f_{V,p}(\gamma_{v}(t))}{n-1}} \left(\bigl(e^{\frac{2(1-\eps)f_{V,p}}{n-1}}\,\Delta_{V}d_{p}\bigl)(\gamma_{v}(t))   \right)^{2}.
\end{align}
\end{lem}
\begin{pf}
When $N=n$,
the vector field $V$ vanishes,
and the desired inequality \eqref{eq:Riccati} is well-known.
In the case of $N=1$ with $\eps =0$,
\eqref{eq:Riccati} has been proved in \cite{KS} (see \cite[Lemma 3.1]{KS}).
We may assume $N\neq 1,n$.

Define $h_{V,v}:=\left(\Delta_{V}d_{p}\right) \circ \gamma_{v}$ and $f_{v}:=f_{V,p} \circ \gamma_{v}$.
By applying the well-known (weighted) Bochner formula to the distance function $d_{p}$,
and by the Cauchy-Schwarz inequality,
\begin{align}\notag
0 &= \Ric^{\infty}_{V}(\dot{\gamma}_{v}(t))+\Vert \nabla^{2} d_{p} \Vert^{2}\left(\gamma_{v}(t)\right)+g(\nabla \Delta_{V}d_p,\nabla d_p)\\  \notag
&\geq \Ric^{N}_{V}(\dot{\gamma}_{v}(t))+\frac{f'_{v}(t)^{2}}{N-n}+\frac{\left(h_{V,v}(t)+f'_{v}(t)\right)^{2}}{n-1}+h'_{V,v}(t)\\  \notag
   &   =   \Ric^{N}_{V}(\dot{\gamma}_{v}(t))+ c\,h^{2}_{V,v}(t)+e^{-\frac{2(1-\eps)f_{v}(t)}{n-1}}\,\left(   e^{\frac{2(1-\eps)f_{v}(t)}{n-1}}  \,h_{V,v}(t)\right)' \\ \notag
   &\qquad \,\,+\frac{1}{n-1}\left( \sqrt{\frac{N-1}{N-n}}f'_{v}(t) +\eps \sqrt{\frac{N-n}{N-1}}h_{V,v}(t)  \right)^2\\ \label{eq:Bochner}
   &\geq \Ric^{N}_{V}(\dot{\gamma}_{v}(t))+ c\,h^{2}_{V,v}(t)+e^{-\frac{2(1-\eps)f_{v}(t)}{n-1}}\,\left(   e^{\frac{2(1-\eps)f_{v}(t)}{n-1}}  \,h_{V,v}(t)\right)'.
\end{align}
We arrive at the desired inequality \eqref{eq:Riccati}.
\end{pf}

\begin{remark}\label{rem:equalRiccati}
{\rm When $N\neq 1,n$,
we assume that
the equality in \eqref{eq:Riccati} holds at $t_{0}\in ]0,\tau(v)[$.
Then we see
\begin{equation*}
\sqrt{\frac{N-1}{N-n}}f'_{v}(t_0) +\eps \sqrt{\frac{N-n}{N-1}}h_{V,v}(t_0)=0
\end{equation*}
since the equality holds in \eqref{eq:Bochner}.
In particular,
if $\eps =0$, then $f'_v(t_0)=0$,
and if $\eps \neq 0$,
\begin{equation*}
h_{V,v}(t_0)=-\eps^{-1}\frac{N-1}{N-n}f'_v(t_0).
\end{equation*}}
\end{remark}

%%%%%%%%%%%%%%%%%%%%%
\subsection{Laplacian comparison theorem}
We now present a Laplacian comparison theorem.
We denote by $\s_{\kappa}(s)$ a unique solution to the Jacobi equation $\psi''(s)+\kappa(s)\,\psi(s)=0$ with $\psi(0)=0$ and $\psi'(0)=1$.
We set
\begin{equation*}
C_{\kappa}:=\inf\{s>0\mid \s_{\kappa}(s)=0\},\quad \cot_{\kappa}(s):=\frac{\s'_{\kappa}(s)}{\s_{\kappa}(s)},\quad H_{\kappa}(s):=c^{-1}\,\cot_{\kappa}(s).
\end{equation*}
Note that
$\cot_{\kappa}(s)$ is a unique solution to the following Riccati equation:
\begin{equation}\label{eq:model Riccati}
\psi'(s)=-\kappa(s)-\psi(s)^2,\quad \lim_{s\downarrow 0}s\psi(s)=1,\quad \lim_{s\uparrow C_{\kappa}}(s-C_{\kappa})\, \psi(s)=1
\end{equation}
under $C_{\kappa}<+\infty$.
We also notice that
$H_{\kappa}$ is decreasing in the case where $\kappa$ is non-negative,
and $H_{\kappa}$ is strictly decreasing in the case where $\kappa$ is a positive function,
or $\kappa$ is a constant function (cf. \cite[Lemma 7.1]{KS}).

%For $\delta \in \mathbb{R}$,
%we also notice the following scaling properties:
%$(1)$ $C_{\kappa\,e^{-4\delta}}=e^{2\delta}\,C_{\kappa}$;
%$(2)$ $\s_{\kappa\,e^{-4\delta}}(s)=e^{2\delta}\,\s_{\kappa}(e^{-2\delta}\,s)$ on $(0,C_{\kappa\,e^{-4\delta}})$;
%$(3)$ $\cot_{\kappa\,e^{-4\delta}}(s)=e^{-2\delta}\,\cot_{\kappa}(e^{-2\delta}\,s)$ on $(0,C_{\kappa\,e^{-4\delta}})$.
%If $\kappa$ is constant,
%then we see $\cot'_{\kappa}(s)<0$ on $(0,\delta_{\kappa})$.

Let us define functions $s_{V,v}:[0,+\infty]\to [0,s_{V,v}(+\infty)]$ and $\tau_V:U_pM\to ]0,+\infty]$ by
\begin{equation}\label{eq:repara}
s_{V,v}(t):=c_p\int^t_0\,e^{      -\frac{2(1-\eps)f_{V,p}(\gamma_v(\xi))}{n-1}       }\,\d\xi,  \quad \tau_{V}(v):=s_{V,v}(\tau(v)).
\end{equation}
Let $t_{V,v}:[0,s_{V,v}(+\infty)]\to [0,+\infty]$ stand for the inverse function of $s_{V,v}$.
We now derive our Laplacian comparison from the Riccati inequality \eqref{eq:Riccati}, which has been obtained by Lu-Minguzzi-Ohta \cite{LMO:CompaFinsler} via Bishop inequality in the gradient case of $V=\nabla f$ (see \cite[Proposition 3.5, Remark 3.10]{LMO:CompaFinsler}):
\begin{thm}\label{lem:Laplacian comparison}
Assume
$\Ric^{N}_{V}(\dot{\gamma}_{v}(t))\geq c^{-1}c^2_p\,\kappa(s_{V,v}(t))\,e^{   -\frac{4(1-\eps)f_{V,p}(\gamma_v(t))}{n-1}      }$ for all $t \in ]0,\tau(v)[$.
Then for all $t \in ]0,\tau(v)[$ with $s_{V,v}(t) \in ]0,\min\{\tau_{V}(v),C_{\kappa} \}[$,
we have
\begin{equation}\label{eq:Laplacian comparison}
\Delta_{V}d_{p}(\gamma_{v}(t)) \leq c_p \,H_{\kappa}(s_{V,v}(t))\,e^{      -\frac{2(1-\eps)f_{V,p}(\gamma_v(t))}{n-1}       }.
\end{equation}
\end{thm}
\begin{pf}
We define two functions $F_{v}:]0,\tau(v)[ \to \mathbb{R}$ and $\hat{F}_{v}:]0,\tau_{V}(v)[ \to \mathbb{R}$ by 
\begin{equation*}
F_{v}:=\bigl(c^{-1}_p\,e^{\frac{2(1-\eps)f_{V,p}}{n-1}}\,\Delta_{V}d_{p}\bigl) \circ \gamma_{v},\quad \hat{F}_{v}:=F_{v}\circ t_{V,v}.
\end{equation*}
From \eqref{eq:Riccati} and the curvature assumption,
for all $s \in ]0,\tau_{V}(v)[$,
\begin{align*}\label{eq:use of Riccati}
\hat{F}'_{v}(s)&    =   F'_{v}(t_{V,v}(s))\, e^{\frac{2(1-\eps)f_{V,p}\left(\gamma_{v}\left(    t_{V,v}(s)  \right)\right) }{n-1}}
c_p^{-1}
\\ \notag
                      & \leq -\Ric^{N}_{V}(\dot{\gamma}_{v}(t_{V,v}(s)))\,e^{\frac{4(1-\eps)f_{V,p}\left(\gamma_{v}\left(    t_{V,v}(s)  \right)\right) }{n-1}}
                      c_p^{-2}
                      -c\,F^{2}_{v}(t_{V,v}(s))\\\notag
                      & \leq -c^{-1}\kappa(s)-c\,\hat{F}^{2}_{v}(s).
\end{align*}
The Riccati equation (\ref{eq:model Riccati}) implies that
for all $s \in ]0,\min\{\tau_{V}(v) ,C_{\kappa} \}  [$
\begin{equation}\label{eq:sharp Riccati}
\hat{F}'_{v}(s)-H'_{\kappa}(s)\leq -c\left(\hat{F}^{2}_{v}(s)-H^2_{\kappa}(s) \right).
\end{equation}

Let us consider a function $G_{\kappa,v}:]0,\min\{\tau_{V}(v) ,C_{\kappa} \}  [\to \mathbb{R}$ by
\begin{equation*}
G_{\kappa,v}:=\mathfrak{s}^{2}_{\kappa}\bigl( \hat{F}_{v}-H_{\kappa} \bigl).
\end{equation*}
From (\ref{eq:sharp Riccati})
it follows that
\begin{align*}\label{eq:monotonicity}
G'_{\kappa,v} &   =  2\,\mathfrak{s}_{\kappa}\, \mathfrak{s}'_{\kappa}\bigl( \hat{F}_{v}-H_{\kappa} \bigl)+\mathfrak{s}^{2}_{\kappa}\bigl( \hat{F}'_{v}-H'_{\kappa} \bigl)\\ \notag
                                    &\leq 2\,\mathfrak{s}_{\kappa}\, \mathfrak{s}'_{\kappa}\bigl(\hat{F}_{v}-H_{\kappa} \bigl)-c\,\mathfrak{s}^{2}_{\kappa}\left(\hat{F}^{2}_{v}-H^{2}_{\kappa}\right)\\ \notag
                                    &  =   -c\,\mathfrak{s}^{2}_{\kappa}\bigl(\hat{F}_{v}-H_{\kappa} \bigl)^{2}\leq 0.   
\end{align*}
Since we see $G_{\kappa,v}(s)\to 0$ as $s\to 0$ by (\ref{eq:model Riccati}),
the function $G_{\kappa,v}$ is non-positive;
in particular,
$\hat{F}_{v} \leq H_{\kappa}$ holds on $]0,\min\{\tau_{V}(v) ,C_{\kappa} \}  [$.
This proves (\ref{eq:Laplacian comparison}).
\end{pf}

\begin{remark}\label{rem:equalLaplacian}
{\rm We assume that
the equality in $(\ref{eq:Laplacian comparison})$ holds at $t_{0}$.
Then $G_{\kappa,v}(s_{V,v}(t_0))=0$.
From $G'_{\kappa,v} \leq 0$
it follows that $G_{\kappa,v}=0$ on $]0,s_{V,v}(t_0)]$;
in particular,
the equality in (\ref{eq:Laplacian comparison}) holds on $]0,t_0]$.}
\end{remark}

\begin{remark}\label{rem:refLaplacian}
{\rm In the gradient case of $V=\nabla f$, this Laplacian comparison theorem has been obtained by Wylie-Yeroshkin \cite{WyYero} under \eqref{eq:RicciLowerBddWY},
the first named author and Li \cite{KL} under \eqref{eq:RicciLowerBddKL},
and Lu-Minguzzi-Ohta \cite{LMO:CompaFinsler} under \eqref{eq:LMOLowerBddd} (see \cite[Theorem 4.4]{WyYero}, \cite[Theorem 2.4]{KL}, and \cite[Remark 3.10]{LMO:CompaFinsler}).
In the non-gradient case, it has been done by Bakry-Qian \cite{BQ2} in the case of $N\in [n,+\infty[$ and $\eps=1$,
and the first named author and Shukuri \cite{KS} in the case of $N\in ]-\infty,1]$ and $\eps=\eps_0$ (see \cite[Theorem 4.2]{BQ2}, and \cite[Theorem 2.5]{KS}).}
\end{remark}

Theorem \ref{lem:Laplacian comparison} leads us to the following:
\begin{lem}\label{lem:Cut point comparisons}
Assume $C_{\kappa}<+\infty$,
and $\Ric^{N}_{V}(\dot{\gamma}_{v}(t))\geq c^{-1}c^2_p\kappa(s_{V,v}(t))\,e^{-\frac{4(1-\eps)f_{V,p}(\gamma_{v}(t))}{n-1}}$ for all $t \in ]0,\tau(v)[$.
Then we have
\begin{equation}\label{eq:Cut point comparisons}
\tau_{V}(v) \leq C_{\kappa}.
\end{equation}
\end{lem}
\begin{pf}
The proof is by contradiction.
Assume $\tau_{V}(v)>C_{\kappa}$.
In this case,
$\tau(v)>t_{V,v}(C_{\kappa})$.
In virtue of (\ref{eq:Laplacian comparison}),
\begin{equation*}
\Delta_{V}d_p(\gamma_{v}(t)) \leq c_p\,H_{\kappa}(s_{V,v}(t))\, e^{-\frac{2(1-\eps)f_{V,p}(\gamma_{v}(t))}{n-1}}
\end{equation*}
for every $t \in ]0,t_{V,v}(C_{\kappa})[$;
in particular,
$\Delta_{V}d_p(\gamma_{v}(t)) \to -\infty$ as $t\to t_{V,v}(C_{\kappa})$ by (\ref{eq:model Riccati}).
This contradicts with the smoothness of $d_{p} \circ \gamma_{v}$ on $]0,\tau(v)[$,
and hence (\ref{eq:Cut point comparisons}).
\end{pf}

\begin{remark}
{\rm Due to Lemma \ref{lem:Cut point comparisons},
one can drop the restriction $s_{V,v}(t) \in ]0,\min\{\tau_{V}(v),C_{\kappa} \}[$ in Theorem \ref{lem:Laplacian comparison}.}
\end{remark}

%%%%%%%%%%%%%%%%%%%%%
\subsection{Rigidity of Laplacian comparison}

We next investigate the equality case of the Laplacian comparison theorem.
\begin{lem}\label{lem:LaplacianRigidity}
Under the same setting as in Theorem \ref{lem:Laplacian comparison},
assume that
the equality in $(\ref{eq:Laplacian comparison})$ holds at $t_{0}\in ]0, \tau(v)  [$.
Choose an orthonormal basis $\{e_{v,i}\}_{i=1}^{n}$ of $T_{p}M$ with $e_{v,n}=v$.
Let $\{Y_{v,i}\}^{n-1}_{i=1}$ and $\{E_{v,i}\}^{n-1}_{i=1}$ be the Jacobi fields and parallel vector fields along $\gamma_{v}$ with $Y_{v,i}(0)=0_p,\,Y_{v,i}'(0)=e_{v,i}$ and $E_{v,i}(0)=e_{v,i}$,
respectively.
Then the following properties hold on $[0,t_0]$:
\begin{enumerate}\setlength{\itemsep}{+0.7mm}
\item If $N=n$, then 
\begin{equation*}
\quad V\equiv 0,\quad Y_{v,i}(t)=\s_{c^2_p\kappa}(t)E_{v,i}(t);
\end{equation*} \label{enum:twisted curv}
\item if $N=1$, then
\begin{equation*}
\eps=0,\quad Y_{v,i}(t)=c^{-1}_p\,\exp\left( \frac{f_{V,p}(\gamma_{v}(t))}{n-1} \right)\,\s_{\kappa}(s_{V,v}(t))\,E_{v,i}(t);
\end{equation*} \label{enum:curv cond}
\item if $N\neq 1,n$, then
\begin{equation*}
\eps=0,\quad g(V,\dot{\gamma}_v(t))\equiv 0,\quad Y_{v,i}(t)=\s_{c^2_p\,\kappa}(t)E_{v,i}(t).
\end{equation*}\label{enum:relaxed twisted curv}
\end{enumerate}
\end{lem}
\begin{pf}
If $N=n$,
then $V\equiv 0$ by definition and the rigidity 
 of Jacobi fields under $\eps=1$ is well-known. For general $\eps\in\R$, its proof can be similarly done (see the proof for (iii) below).
If $N=1$,
then the desired assertion has been proved by Kuwae-Shukuri \cite{KS} (see \cite[Lemma 3.2]{KS}, and also \cite[Lemma 4.3]{{WyYero}}).
We may assume $N\neq 1,n$.

We first show $\eps = 0$ by contradiction.
We suppose $\eps \neq 0$.
Then in view of Remarks \ref{rem:equalRiccati} and \ref{rem:equalLaplacian},
$h_{V,v}(t)$ is equal to
\begin{equation*}
-\eps^{-1}\frac{N-1}{N-n}f'_v(t)=c_p \,H_{\kappa}(s_{V,v}(t))\,e^{      -\frac{2(1-\eps)f_v(t)}{n-1}       }
\end{equation*}
on $]0,t_0]$,
where $h_{V,v}:=\left(\Delta_{V}d_{p}\right) \circ \gamma_{v}$ and $f_{v}:=f_{V,p} \circ \gamma_{v}$.
This is a contradiction since the left hand side converges as $t\to 0$,
but the right hand side does not.
We now possess $\eps =0$.
Remark \ref{rem:equalRiccati} says $g(V,\dot{\gamma}_v(t))\equiv 0$;
moreover,
\begin{align*}
\Ric^{1}_{V}(\dot{\gamma}_{v}(t))&\geq \Ric^{N}_{V}(\dot{\gamma}_{v}(t))\geq (n-1)c^2_p\kappa(s_{V,v}(t))\,e^{-\frac{4f_{V,p}(\gamma_{v}(t))}{n-1}},\\
\Delta_{V}d_{p}(\gamma_{v}(t))&=c_p \,H_{\kappa}(s_{V,v}(t))\,e^{      -\frac{2f_v(t)}{n-1}       },
\end{align*}
and hence the equality for $N=1$ occurs.
By the rigidity of Jacobi fields for $N=1$,
and $g(V,\dot{\gamma}_v(t))\equiv 0$,
we conclude
\begin{equation*}
Y_{v,i}(t)=c^{-1}_p\,\s_{\kappa}(c_{p}t)\,E_{v,i}(t)=\s_{c^2_p\, \kappa}(t)E_{v,i}(t).
\end{equation*}
Thus we complete the proof.
\end{pf}

%%%%%%%%%%%%%%%%%%%%%%%%
\subsection{Laplacian comparison with bounded density}\label{sec:Bounded densities}
In this last subsection,
we focus on the gradient case of $V=\nabla f$.
We show comparison results under a bound for $f$.
As stated above,
we choose $c_p$ as \eqref{eq:gradcp}.

\begin{lem}\label{lem:finite cut value comparison}
Let $V=\nabla f$.
We assume
\begin{equation*}
C_{\kappa}<+\infty,\quad \Ric^{N}_{f}(\dot{\gamma}_{v}(t))\geq c^{-1}\kappa(s_{\nabla f,v}(t))\,e^{-\frac{4(1-\eps)f(\gamma_{v}(t))}{n-1}},\quad
(1-\eps)f\circ \gamma_{v} \leq (n-1)\delta
\end{equation*}
on $]0,\tau(v)[$ for $\delta \in \mathbb{R}$.
Then
\begin{equation*}\label{eq:finite Cut point comparisons}
\tau(v)\leq C_{ \kappa e^{-4\delta}}.
\end{equation*}
\end{lem}
\begin{pf}
The upper bound for $f$ implies $e^{-2\delta}\,\tau(v) \leq \tau_{\nabla f}(v)$.
By Lemma \ref{lem:Cut point comparisons} and $C_{\kappa e^{-4\delta}}=e^{2\delta}\,C_{\kappa}$,
we complete the proof.
\end{pf}

\begin{remark}\label{rem:refLaplacian}
{\rm In the gradient case of $V=\nabla f$, a similar upper bound has been shown by Lu-Minguzzi-Ohta \cite{LMO:CompaFinsler} (see \cite[Theorem 3.6]{LMO:CompaFinsler}).}
\end{remark}

We write $\Delta_f:=\Delta_{\nabla f}$.
In view of Lemma \ref{lem:finite cut value comparison},
we have the following:
\begin{lem}\label{lem:finite Laplacian comparison}
Let $V=\nabla f$.
We assume
\begin{equation*}
\Ric^{N}_{f}(\dot{\gamma}_{v}(t))\geq c^{-1}\kappa(s_{\nabla f,v}(t))\,e^{-\frac{4(1-\eps)f(\gamma_{v}(t))}{n-1}},\quad
(1-\eps)f\circ \gamma_{v} \leq (n-1)\delta
\end{equation*}
on $]0,\tau(v)[$ for $\delta \in \mathbb{R}$.
We further assume that
$H_{\kappa}$ is decreasing.
Then for all $t \in ]0,\tau(v)[$,
\begin{equation}\label{eq:finite Laplacian comparison}
\Delta_{f}\, d_{p}(\gamma_{v}(t))\leq e^{2\delta}\,H_{\kappa\,e^{-4\delta}}(t)\,e^{-\frac{2(1-\eps)f(\gamma_{v}(t))}{n-1}}.
\end{equation}
\end{lem}
\begin{pf}
From the upper boundedness of $f$,
we deduce $s_{\nabla f,v}(t)\geq e^{-2\delta}t$ for every $t \in ]0,\tau(v)[$.
Now,
the assumption for $H_{\kappa}$ and (\ref{eq:Laplacian comparison}) imply
\begin{equation}\label{eq:proof of finite Laplacian comparison}
\Delta_{f}\, d_{p}(\gamma_{v}(t))\leq H_{\kappa}(s_{\nabla f,v}(t))\,e^{-\frac{2(1-\eps)f(\gamma_{v}(t))}{n-1}} \leq H_{\kappa}\left(e^{-2\delta}t\right)\,e^{-\frac{2(1-\eps)f(\gamma_{v}(t))}{n-1}}.
\end{equation}
The right hand side is equal to that of (\ref{eq:finite Laplacian comparison}).
\end{pf}

\begin{remark}\label{rem:equality in finite Laplacian comparison}
{\rm Assume that
the equality in $(\ref{eq:finite Laplacian comparison})$ holds at $t_{0}\in ]0,\tau(v)[$.
Then the equalities in (\ref{eq:proof of finite Laplacian comparison}) also hold.
The equality in (\ref{eq:Laplacian comparison}) holds (see Lemma \ref{lem:LaplacianRigidity}).
Moreover, if $H_{\kappa}$ is strictly decreasing, then
$s_{\nabla f,v}(t_{0})=e^{-2\delta}t_{0}$ 
%by the assumption for $H_{\kappa}$,
and hence $(1-\eps)f\circ \gamma_{v} = (n-1)\delta$ on $[0,t_{0}]$.}
\end{remark}

\begin{remark}\label{rem:refLaplacian}
{\rm In the gradient case of $V=\nabla f$, a similar estimate has been shown by Lu-Minguzzi-Ohta \cite{LMO:CompaFinsler} (see \cite[Theorem 3.9]{LMO:CompaFinsler}).
Here they further assumed a lower bound of $(1-\eps)f\circ \gamma_{v}$, and obtained an estimate that does not depend on $f$.
Moreover,
they have concluded a volume estimate under the same setting (see \cite[Theorem 3.11]{LMO:CompaFinsler}).}
\end{remark}

%%%%%%%%%%%%%%%%%%%%%
%%%%%%%%%%%%%%%%%%%%%
%%%%%%%%%%%%%%%%%%%%%
\section{Diameter}\label{sec:Diameter}

%%%%%%%%%%%%%%%%%%%%%
\subsection{Diameter comparison theorem}

We present the following comparison of Bonnet-Myers type:
\begin{prop}\label{thm:diameter comparison}
We assume $C_{\kappa}<+\infty$,
and also assume $\Ric^{N}_{V} \geq c^{-1}c^2_p\,\kappa_{V,p}\,e^{   -\frac{4(1-\eps)f_{V,p}}{n-1}      }g$.
Then we have
\begin{equation*}\label{eq:diameter comparison}
\sup_{x\in M} s_{V,p}(x) \leq C_{\kappa}.
\end{equation*}
\end{prop}
\begin{pf}
Fix $x\in M$,
and a unit speed minimal geodesic $\gamma:[0,d(p,x)] \to M$ from $p$ to $x$.
According to Lemma \ref{lem:Cut point comparisons},
we see
\begin{equation*}
s_{V,p}(x) \leq c_p\int^{d_p(x)}_{0}\, e^{      -\frac{2(1-\eps)f_{V,p}(\gamma_v(\xi))}{n-1}} \,\d\xi \leq \tau_{V}(v)\leq C_{\kappa},
\end{equation*}
where $v:=\dot{\gamma}(0)$.
We conclude the desired assertion.
\end{pf}

By the same method,
we can deduce the following from Lemma \ref{lem:finite cut value comparison}:
\begin{prop}\label{thm:finite diameter comparison}
Let $V=\nabla f$.
We assume
\begin{equation*}
C_{\kappa}<+\infty,\quad \Ric^{N}_{f} \geq c^{-1}\,\kappa_{\nabla f,p}\,e^{   -\frac{4(1-\eps)f}{n-1}      }g,\quad (1-\eps)f \leq (n-1)\delta
\end{equation*}
for $\delta \in \R$.
Then
\begin{equation*}\label{eq:finite diameter comparison}
\sup_{x\in M} d_p(x) \leq C_{\kappa e^{-4\delta}}.
\end{equation*}
In particular,
$M$ is compact.
\end{prop}

%%%%%%%%%%%%%%%%%%%%%
\subsection{Maximal diameter theorem}
We now establish a maximal diameter theorem for the equality case of Proposition \ref{thm:diameter comparison}.
We consider a conformally deformed Riemannian metric (with singularity) defined by
\begin{equation*}
g_{V,p}:=c^2_p\,e^{-\frac{4(1-\eps)f_{V,p}}{n-1}} g,
\end{equation*}
and its distance function
\begin{align*}
d_{g_{V,p}}(x,y)=c_p\inf_{\sigma}\int^{l}_{0}\,e^{-\frac{2(1-\eps)f_{V,p}(\sigma(\xi))}{n-1}}g(\dot{\sigma}(\xi),\dot{\sigma}(\xi))^{1/2}\,\d\xi,
\end{align*}
where the infimum is taken over all piecewise smooth curves $\sigma:[0,l]\to M$ with $\sigma(0)=x$ and $\sigma(l)=y$.
Note that this satisfies the triangle inequality,
and $d_{g_{V,p}}(p,x)\leq s_{V,p}(x)$ for all $x\in M$.
The following is one of our main theorems:
\begin{thm}\label{thm:diameter rigidity}
We assume $C_{\kappa}<+\infty$,
and also assume $\Ric^{N}_{V} \geq c^{-1}c^2_p\,\kappa_{V,p}\,e^{   -\frac{4(1-\eps)f_{V,p}}{n-1}      }g$.
Then we have
\begin{equation}\label{eq:diameter rigidity comparison}
\sup_{x\in M} d_{g_{V,p}}(p,x) \leq C_{\kappa}.
\end{equation}
Moreover,
we further assume that $\kappa(s)=\kappa(C_{\kappa}-s)$ for all $s\in [0,C_{\kappa}]$,
and $\kappa$ is positive.
If there exists $q \in M$ with
\begin{align}\label{eq:maxdiamas}
\Ric^{N}_{V} \geq c^{-1}c^2_q\,\kappa_{V,q}\,e^{   -\frac{4(1-\eps)f_{V,q}}{n-1}      }g,\quad c_q\,e^{   -\frac{2(1-\eps)f_{V,q}}{n-1}      }=c_p\,e^{   -\frac{2(1-\eps)f_{V,p}}{n-1}      }
\end{align}
such that
\begin{equation*}
d_{g_{V,p}}(p,q) = C_{\kappa},
\end{equation*}
then
\begin{equation*}
d_p+d_q \equiv d(p,q)
\end{equation*}
on $M$,
and by identifying $U_{p}M$ with the $(n-1)$-dimensional unit sphere $(\mathbb{S}^{n-1},g_{\mathbb{S}^{n-1}})$,
we have the following rigidity properties:
\begin{enumerate}\setlength{\itemsep}{+1.0mm}
\item If $N=n$, then $V\equiv 0$, and $g=dt^2+\s_{c^2_p\kappa}(t)\,g_{\mathbb{S}^{n-1}}$;
\item if $N=1$, then $\eps=0$, and
\begin{equation*}
g=dt^2+c^{-2}_p\,\exp\left( \frac{2f_{V,p}(\gamma_{v}(t))}{n-1} \right)\,\s^2_{\kappa}(s_{V,v}(t))g_{\mathbb{S}^{n-1}};
\end{equation*}
\item if $N\neq 1,n$, then $\eps=0$, $V$ is orthogonal to $\nabla d_p$ on $M\setminus \{p,q\}$ and vanishes at $\{p,q\}$, and $g=dt^2+\s_{c^2_p\kappa}(t)\,g_{\mathbb{S}^{n-1}}$.
\end{enumerate}
\end{thm}
\begin{pf}
The inequality (\ref{eq:diameter rigidity comparison}) is a direct consequence of Proposition \ref{thm:diameter comparison}.
Let us prove the rigidity part.
Set
\begin{equation}\label{eq:triangle inequality holding domain}
\Omega_{p,q}:=\left\{ \,x \in M \setminus \{p,q\}\, \middle| \,d_{p}(x)+d_{q}(x)=d(p,q)\, \right\}.
\end{equation}
The interior of a unit speed minimal geodesic from $p$ to $q$ lies in $\Omega_{p,q}$,
and hence
$\Omega_{p,q}$ is a non-empty closed subset of $M \setminus \{p,q\}$.

We show that
$\Omega_{p,q}$ is open.
Fix $x\in \Omega_{p,q}$.
Note that
$x$ does not belong to the cut locus of $p$ and $q$.
We take a sufficiently small domain $\Omega\subset M$ containing $x$ on which  
$d_{p}$ and $d_{q}$ are smooth.
We apply Theorem \ref{lem:Laplacian comparison} to them with the help of the first assumption in \eqref{eq:maxdiamas}.
By using the second one in \eqref{eq:maxdiamas},
for each $y \in \Omega$ we see
\begin{align*}
\Delta_{V}(d_p+d_q)(y)&\leq c_p \,H_{\kappa}(s_{V,p}(y))\,e^{      -\frac{2(1-\eps)f_{V,p}(y)}{n-1}       }+c_q \,H_{\kappa}(s_{V,q}(y))\,e^{      -\frac{2(1-\eps)f_{V,q}(y)}{n-1}       }\\
                                     &= c^{-1}c_p \,\left(\cot_{\kappa}(s_{V,p}(y))+\cot_{\kappa}(s_{V,q}(y))\right)\,e^{      -\frac{2(1-\eps)f_{V,p}(y)}{n-1}       }.
\end{align*}
The second one in \eqref{eq:maxdiamas} implies $g_{V,p}=g_{V,q}$.
The triangle inequality for $d_{g_{V,p}}$ leads to
\begin{align*}\label{eq:weighted triangle inequality}
s_{V,p}(y)+s_{V,q}(y)& \geq d_{g_{V,p}}(p,y)+d_{g_{V,q}}(q,y)= d_{g_{V,p}}(p,y)+d_{g_{V,p}}(q,y)\geq d_{g_{V,p}}(p,q)=C_{\kappa}.
\end{align*}
We now recall that 
if $\kappa(s)=\kappa(\delta_{\kappa}-s)$,
then $\s_{\kappa}(s)=\s_{\kappa}(C_{\kappa}-s)$,
and also
if $\kappa$ is positive,
then $\cot_{\kappa}$ is strictly decreasing (cf. \cite[Lemma 7.1]{KS}).
Therefore,
\begin{equation*}
\cot_{\kappa}(s_{V,q}(y))\leq \cot_{\kappa}(C_{\kappa}-s_{V,p}(y))=-\cot_{\kappa}(s_{V,p}(y)).
\end{equation*}
It follows that $\Delta_{V}(d_p+d_q)(y)\leq 0$.
Applying the strong maximum principle over $\Omega$,
we conclude $\Omega \subset \Omega_{p,q}$ and the openness.

From the connectedness of $M\setminus \{p,q\}$, 
we conclude $\Omega_{p,q}=M\setminus \{p,q\}$.
Moreover, by using the second one in \eqref{eq:maxdiamas} again,  
the equality in (\ref{eq:Laplacian comparison}) holds on $M\setminus \{p,q\}$;
in particular,
Lemma \ref{lem:LaplacianRigidity} yields the rigidity properties of $\eps,V,g$ via a homeomorphism $\Theta:[0,d(p,q)]\times U_{p}M\to M$ defined by $\Theta(t,v):=\gamma_{v}(t)$.
Thus,
we complete the proof.
\end{pf}

In the gradient case of $V=\nabla f$,
the metric $g_{\nabla f,p}$ can be written as $g_f$.
In that case,
and the case where $\kappa$ is constant,
we conclude:
\begin{cor}\label{cor:corodiameter rigidity}
Let us assume $V=\nabla f$,
and let $\kappa$ be a positive constant.
We assume $\Ric^{N}_{V} \geq c^{-1}\kappa\,e^{   -\frac{4(1-\eps)f}{n-1}      }g$.
Then we have
\begin{equation*}
\sup_{x\in M} d_{g_{f}}(p,x) \leq C_{\kappa}.
\end{equation*}
Moreover,
if there exists $q \in M$ such that
\begin{equation*}
d_{g_{f}}(p,q) = C_{\kappa},
\end{equation*}
then the following rigidity properties hold:
\begin{enumerate}\setlength{\itemsep}{+1.0mm}
\item If $N=n$, then $f$ is constant, and $M$ is isometric to a sphere with constant curvature $\kappa\,e^{   -\frac{4(1-\eps)f}{n-1}      }$;
\item if $N=1$, then $\eps=0$, $f$ is radial with respect to $p$ $($i.e., $f$ depends only on $d_p)$, $M$ is homeomorphic to a sphere, and
\begin{equation*}
g=dt^2+\exp\left( 2\frac{f(\gamma_{v}(t))+f(p)}{n-1} \right)\,\s^2_{\kappa}\left(s_{\nabla f,v}(t)\right)g_{\mathbb{S}^{n-1}};
\end{equation*}
\item if $N\neq 1,n$, then $\eps=0$, $f$ is constant, and $M$ is isometric to a sphere with constant curvature $\kappa\,e^{   -\frac{4 f}{n-1}      }$.
\end{enumerate}
\end{cor}
\begin{pf}
Almost all parts are the direct consequence of Theorem \ref{thm:diameter rigidity}.
Actually,
the assumption (\ref{eq:maxdiamas}) is always satisfied in this setting.
We only need to verify the radial property of $f$ in the case of $N=1$.
This has been proved by Wylie-Yeroshkin (see \cite[Proposition 4.15]{WyYero}).
Thus we arrive at the desired conclusion.
\end{pf}

\begin{remark}\label{rem:refmaxdiam}
{\rm In the gradient case of $V=\nabla f$, Corollary \ref{cor:corodiameter rigidity} has been obtained by Wylie-Yeroshkin \cite{WyYero} under \eqref{eq:RicciLowerBddWY},
and the first named author and Shukuri \cite{KS} under \eqref{eq:RicciLowerBddKL} (see \cite[Theorem 4.16]{WyYero}, and \cite[Corollary 2.22]{KS}).
In the non-gradient case,
Kuwada \cite{K} has proven Theorem \ref{thm:diameter rigidity} in the case where $N\in [n,+\infty[,\,\eps=1$ and $\kappa$ is constant.}
\end{remark}

%%%%%%%%%%%%%%%%%%%%%
\subsection{Maximal diameter theorem with bounded density}

We next investigate rigidity phenomena for the equality case of Proposition \ref{thm:finite diameter comparison},
which is new even in the setting of Wylie-Yeroshkin \cite{WyYero}.
\begin{thm}\label{thm:finite diameter rigidity}
Let $V=\nabla f$.
We assume
\begin{equation*}
C_{\kappa}<+\infty,\quad \Ric^{N}_{f} \geq c^{-1}\,\kappa_{\nabla f,p}\,e^{   -\frac{4(1-\eps)}{n-1}      },\quad (1-\eps)f \leq (n-1)\delta
\end{equation*}
for $\delta \in \R$.
Then
\begin{equation}\label{eq:hogehoge}
\sup_{x\in M} d_p(x) \leq C_{\kappa e^{-4\delta}}.
\end{equation}
Moreover,
we further assume that $\kappa(s)=\kappa(C_{\kappa}-s)$ for all $s\in [0,C_{\kappa}]$,
and $\kappa$ is positive.
If there exists $q\in M$ with
\begin{equation*}\label{eq:finmaxdiamas}
\Ric^{N}_{f} \geq c^{-1}\,\kappa_{\nabla f,q}\,e^{   -\frac{4(1-\eps)f}{n-1}      }g
\end{equation*}
such that
\begin{equation*}
d_p(q)=C_{\kappa e^{-4\delta}},
\end{equation*}
then
\begin{equation}\label{eq:sphere}
d_p+d_q \equiv C_{\kappa e^{-4\delta}}
\end{equation}
on $M$,
and the following rigidity properties hold:
\begin{enumerate}\setlength{\itemsep}{+1.0mm}
\item If $N=n$, then $(1-\eps)f\equiv (n-1)\delta$, and $g=dt^2+\s_{\kappa e^{-4\delta}}(t)\,g_{\mathbb{S}^{n-1}}$;
\item if $N\neq n$, then $\eps=0$, $f\equiv (n-1)\delta$, and $g=dt^2+\s_{\kappa e^{-4\delta}}(t)\,g_{\mathbb{S}^{n-1}}$.
\end{enumerate}
\end{thm}
\begin{pf}
The inequality (\ref{eq:hogehoge}) was proved in Proposition \ref{thm:finite diameter comparison}.
We will prove the rigidity part by using Lemma \ref{lem:finite Laplacian comparison} instead of Theorem \ref{lem:Laplacian comparison} along the lines of the proof of Theorem \ref{thm:diameter rigidity}.
Define a non-empty closed subset $\Omega_{p,q}$ of $M \setminus \{p,q\}$ as (\ref{eq:triangle inequality holding domain}).
We show the openness of $\Omega_{p,q}$.
For a fixed $x\in \Omega_{p,q}$,
take a domain $\Omega\subset M$ containing $x$ on which $d_p$ and $d_q$ are smooth.
Due to Lemma \ref{lem:finite Laplacian comparison} for each $y \in \Omega$,
\begin{align*}
\Delta_{f}(d_{p}+d_{q})(y)&\leq e^{2\delta}\,\left(H_{\kappa\,e^{-4\delta}}(d_p(x))\,e^{-\frac{2(1-\eps)f(y)}{n-1}}+  H_{\kappa\,e^{-4\delta}}(d_{q}(y))\,e^{-\frac{2(1-\eps)f(y)}{n-1}} \right)\\
                                         &= c^{-1}\,e^{2\delta}\,\left(\cot_{\kappa\,e^{-4\delta}}(d_p(x))+  \cot_{\kappa\,e^{-4\delta}}(d_{q}(y)) \right)e^{-\frac{2(1-\eps)f(y)}{n-1}},
\end{align*}
where we notice that
$H_{\kappa}$ is strictly decreasing by the positivity of $\kappa$.
Since it holds that $\s_{\kappa}(s)=\s_{\kappa}(C_{\kappa}-s)$,
we have
\begin{equation*}
\cot_{\kappa\,e^{-4\delta}}(d_{q}(y))\leq \cot_{\kappa\,e^{-4\delta}}(C_{\kappa\,e^{-4\delta}}-d_{p}(y))=-\cot_{ \kappa\,e^{-4\delta}}(d_{p}(y)),
\end{equation*}
and obtain $\Delta_{f}(d_p+d_q)(y)\leq 0$.
According to the strong maximum principle,
we arrive at the openness of $\Omega_{p,q}$.

From the same argument as in the proof of Theorem \ref{thm:diameter rigidity},
one can conclude (\ref{eq:sphere}).
Now,
the equality in (\ref{eq:finite Laplacian comparison}) holds on $M\setminus \{p,q\}$ (see Remark \ref{rem:equality in finite Laplacian comparison}).
We can apply Lemma \ref{lem:LaplacianRigidity} to our situation,
and $(1-\eps)f \equiv (n-1)\delta$.
This completes the proof.
\end{pf}

For constant $\kappa$, we have the following:
\begin{cor}\label{cor:cor finite diameter rigidity}
Let $V=\nabla f$,
and let $\kappa$ be a positive constant.
We assume
\begin{equation*}
\Ric^{N}_{f} \geq c^{-1}\,\kappa\,e^{   -\frac{4(1-\eps)f}{n-1}      },\quad (1-\eps)f \leq (n-1)\delta
\end{equation*}
for $\delta \in \R$.
Then
\begin{equation*}\label{eq:finite diameter rigicps}
\sup_{x\in M} d_p(x) \leq C_{\kappa e^{-4\delta}}.
\end{equation*}
Moreover,
If there exists $q\in M$ such that
\begin{equation*}
d_p(q)=C_{\kappa e^{-4\delta}}.
\end{equation*}
then the following rigidity properties hold:
\begin{enumerate}\setlength{\itemsep}{+1.0mm}
\item If $N=n$, then $(1-\eps )f\equiv (n-1)\delta$, and $M$ is isometric to a sphere of constant curvature $\kappa e^{-4\delta}$;
\item if $N\neq n$, then $\eps=0$, $f\equiv (n-1)\delta$, and $M$ is isometric to a sphere of constant curvature $\kappa e^{-4\delta}$.
\end{enumerate}
\end{cor}

%%%%%%%%%%%%%%%%%%%%%
%%%%%%%%%%%%%%%%%%%%%
%%%%%%%%%%%%%%%%%%%%%
\section{Volume}\label{sec:Volume}

%%%%%%%%%%%%%%%%%%%%%%%%%%%%%%
\subsection{Volume elements}\label{sec:volume elements}

For $t \in ]0,\tau(v)[$,
and for the volume element $\theta(t,v)$ of the $t$-level surface of $d_p$ at $\gamma_{v}(t)$,
\begin{equation*}\label{eq:volume element}
\theta_{V}(t,v):=e^{-f_{V,p}(\gamma_{v}(t))}\, \theta(t,v),\quad \hat{\theta}_{V}(s,v):=\theta_{V}(t_{V,v}(s),v),
\end{equation*}
where $t_{V,v}$ is the inverse function of $s_{V,v}$ defined as (\ref{eq:repara}).
We first show:
\begin{lem}\label{lem:volume element comparison}
Assume that
$\Ric^{N}_{V}(\dot{\gamma}_{v}(t))\geq c^{-1}c^2_p\,\kappa(s_{V,v}(t))\,e^{   -\frac{4(1-\eps)f_{V,p}(\gamma_v(t))}{n-1}      }$ for all $t \in ]0,\tau(v)[$.
Then for all $s_{1},s_{2} \in ]0,\tau_{V}(v)[$ with $s_{1}\leq s_{2}$
\begin{equation}\label{eq:hoge}
\frac{\hat{\theta}_{V}(s_{2},v)}{ \hat{\theta}_{V}(s_{1},v)}\leq \frac{\mathfrak{s}^{1/c}_{\kappa}(s_{2})}{\mathfrak{s}^{1/c}_{\kappa}(s_{1})}.
\end{equation}
Moreover,
if $c=(n-1)^{-1}$, then
for all $s\in [0,\tau_{V}(v)[$ we have
\begin{equation}\label{eq:absvolele}
\hat{\theta}_{V}(s,v)\leq \mathfrak{s}^{n-1}_{\kappa}(s).
\end{equation}
\end{lem}
\begin{pf}
Let us use the inequality \eqref{eq:Laplacian comparison}.
For all $s \in ]0,\tau_{V}(v)[$
we see
\begin{equation*}
\frac{\d}{\d s}\log \frac{\hat{\theta}_{V}(s,v)}{\mathfrak{s}^{1/c}_{\kappa}(s)}=\bigl(c^{-1}_p\,e^{\frac{2(1-\eps)f_{V,p}}{n-1}}\,\Delta_{V}d_p\bigl)(\gamma_{v}(t_{V,v}(s)))-H_{\kappa}(s)\leq 0.
\end{equation*}
This implies \eqref{eq:hoge}.
If $c=(n-1)^{-1}$,
then we see $\hat{\theta}_{V}(s,v)/\mathfrak{s}^{1/c}_{\kappa}(s)\to 1$ as $s\to 0$ by \eqref{eq:model Riccati}.
Hence we arrive at \eqref{eq:absvolele}.
\end{pf}

\begin{remark}\label{rem:equality in volume element comparison}
{\rm Assume that
the equality in \eqref{eq:absvolele} holds at $s_{0}\in ]0,\tau_{V}(v)[$.
Then the equality in \eqref{eq:absvolele} holds on $[0,s_{0}]$;
in particular,
the equality in \eqref{eq:Laplacian comparison} holds on $]0,t_{V,v}(s_{0})]$ (see Lemma \ref{lem:LaplacianRigidity}).}
\end{remark}

\begin{remark}
{\rm We have $c=(n-1)^{-1}$ if and only if either (1) $N=n$; or (2) $N=1$; or (3) $N\neq 1,n$ and $\eps=0$ (cf. Lemma \ref{lem:LaplacianRigidity}).}
\end{remark}

%%%%%%%%%%%%%%%%%%%%%%%%%%%%%%
\subsection{Volume comparison theorem}\label{sec:volume comparisons}
For $r>0$,
we define 
\begin{equation*}
B_{V,r}(p):=\left\{\,x\in M \mid s_{V,p}(x) < r \,\right\},
\end{equation*}
and also define measures
\begin{equation*}
\mu_{V,p}:=e^{-f_{V,p}}\,v_g,\quad \nu_{V,p}:=e^{ - \frac{2(1-\eps)f_{V,p}}{n-1}  }\mu_{V,p},
\end{equation*}
where $v_g$ is the Riemannian volume measure.
%\begin{remark}
%{\rm We consider an invariant measure $\mu$ for $\Delta_V$ (see e.g., \cite{BRW:invariant}, \cite{IkedaWatanabe:SDE}), 
%which is a Radon measure on $M$ and is a solution of $\Delta_V^*\mu=0$, namely, for any $\varphi \in C_0^{\infty}(M)$,
%\begin{align*}
%\int_M\Delta_V \varphi\d\mu=0.
%\end{align*}
%One can expect that $\mu_{V,p}$ is an invariant measure for $\Delta_V$ (cf. \cite{BQ2}, \cite[Proposition~4.5]{KS}).}
%\end{remark}
By straightforward argument,
one can verify
\begin{equation}\label{eq:integration formula}
\nu_{V,p}\left(B_{V,r}(p)\right)=\int_{U_{p}M}\,\int^{r}_{0}\,\bar{\theta}_{V}(s,v)\,\d s\, \d v,
\end{equation}
where
\begin{equation*}\label{eq:extended volume element}
\bar{\theta}_{V}(s,v):=\begin{cases}
                                                    \hat{\theta}_{V}(s,v) & \text{if $s<\tau_{V}(v)$}, \\
                                                                      0           & \text{if $s \geq \tau_{V}(v)$}.
                                                   \end{cases}
\end{equation*}
We also set
\begin{equation*}
\mathcal{S}_{\kappa}(r):=\int^{r}_{0}\,\bar{\s}^{1/c}_{\kappa}(s)\,\d s,               
\end{equation*}
where
\begin{equation*}
\bar{\s}_{\kappa}(s):=\begin{cases}
                                                    \s_{\kappa}(s) & \text{if $s<C_{\kappa}$}, \\
                                                                                0           & \text{if $s \geq C_{\kappa}$}.
                                                   \end{cases}                                          
\end{equation*}

We first present the following absolute comparison theorem of Bishop type:
\begin{prop}\label{thm:absolute volume comparison}
We assume $\Ric^{N}_{V} \geq c^{-1}c^2_p\,\kappa_{V,p}\,e^{   -\frac{4(1-\eps)f_{V,p}}{n-1}      }g$,
and $c=(n-1)^{-1}$.
Then for all $r>0$
we have
\begin{equation}\label{eq:absvolcompa}
\nu_{V,p}(B_{V,r}(p))\leq \omega_{n-1}\,\mathcal{S}_{\kappa}(r),
\end{equation}
where $\omega_{n-1}$ is the volume of the $(n-1)$-dimensional unit sphere.
In particular,
\begin{equation*}
\varlimsup_{r\to +\infty}\frac{\nu_{V,p}(B_{V,r}(p))}{\mathcal{S}_{\kappa}(r)}\leq \omega_{n-1}.
\end{equation*}
\end{prop}
\begin{pf}
By \eqref{eq:absvolele} in Lemma \ref{lem:volume element comparison},
for all $s \geq 0$ and $v\in U_{p}M$
\begin{equation}\label{eq:extended absolute volume element comparison}
\bar{\theta}_{V}(s,v) \leq \bar{\s}^{n-1}_{\kappa}(s).
\end{equation}
Integrating it over $]0,r[$ with respect to $s$,
and \eqref{eq:integration formula} complete the proof.
\end{pf}

\begin{remark}\label{rem:equality in absolute volume comparison}
{\rm Assume that
the equality in \eqref{eq:absvolcompa} holds.
Then the equality in (\ref{eq:extended absolute volume element comparison}) holds for all $s\in [0,r]$ and $v\in U_{p}M$.
We have $\tau_{V}(v) \geq \min\{r,C_{\kappa}\}$ for all $v\in U_{p}M$,
and the equality in \eqref{eq:absvolele} holds for all $s\in [0,\min\{r,C_{\kappa}\}[$ and $v\in U_{p}M$ (see Remark~\ref{rem:equality in volume element comparison}).}
\end{remark}

We also prove the following relative comparison of Bishop-Gromov type:
\begin{prop}\label{thm:relative volume comparison}
We assume $\Ric^{N}_{V} \geq c^{-1}c^2_p\,\kappa_{V,p}\,e^{   -\frac{4(1-\eps)f_{V,p}}{n-1}      }g$.
Then for all $r,R>0$ with $r\leq R$
we have
\begin{equation*}
\frac{\nu_{V,p}(B_{V,R}(p))}{\nu_{V,p}(B_{V,r}(p))} \leq \frac{\mathcal{S}_{\kappa}(R)}{\mathcal{S}_{\kappa}(r)}.
\end{equation*}
\end{prop}
\begin{pf}
By using \eqref{eq:hoge},
for all $s_{1},s_{2}> 0$ with $s_{1}\leq s_{2}$,
and $v\in U_{p}M$
\begin{equation*}
\bar{\theta}_{V}(s_{2},v)\; \bar{\s}_{\kappa}^{1/c}(s_{1}) \leq \bar{\theta}_{V}(s_{1},v)\; \bar{\s}_{\kappa}^{1/c}(s_{2}).
\end{equation*}
Let us integrate the both sides over $]0,r[$ with respect to $s_{1}$,
and over $]r,R[$ with respect to $s_{2}$.
We obtain
\begin{equation*}
\frac{\int^{R}_{r}\bar{\theta}_{V}(s_{2},v)\,\d s_{2}}{\int^{r}_{0}\bar{\theta}_{V}(s_{1},v)\,\d s_{1}}\leq \frac{\mathcal{S}_{\kappa}(R)-\mathcal{S}_{\kappa}(r)}{\mathcal{S}_{\kappa}(r)}.
\end{equation*}
The formula \eqref{eq:integration formula} yields
\begin{equation*}
          \frac{\nu_{V,p}( B_{V,R}(p))}{\nu_{V,p}(B_{V,r}(p))}
  =   1+\frac{\int_{U_{p}M} \int^{R}_{r}\,\bar{\theta}_{V}(s_{2},v)\,\d s_{2}\,\d v}{\int_{U_{p}M} \int^{r}_{0}\bar{\theta}_{V}(s_{1},v)\,\d s_{1}\,\d v}\leq \frac{\mathcal{S}_{\kappa}(R)}{\mathcal{S}_{\kappa}(r)}.
\end{equation*}
We complete the proof of Proposition \ref{thm:relative volume comparison}.
\end{pf}

\begin{remark}\label{rem:refLaplacian}
{\rm In the gradient case of $V=\nabla f$, similar volume comparison theorems have been studied by Wylie-Yeroshkin \cite{WyYero} under \eqref{eq:RicciLowerBddWY},
and by the first named author and Li \cite{KL} under \eqref{eq:RicciLowerBddKL} (see \cite[Corollary 4.6]{WyYero}, \cite[Theorem 2.10]{KL}).
In the non-gradient case, it has been done by the first named author and Shukuri \cite{KS} in the case of $N\in ]-\infty,1]$ and $\eps=\eps_0$ (see \cite[Theorem 2.14]{KS}).}
\end{remark}

%%%%%%%%%%%%%%%%%%%%%%%%%%%%
\subsection{Rigidity of volume comparison}\label{sec:Volume growth rigidity}
We investigate the equality cases of volume comparisons (cf. \cite[Theorem 4.17]{WyYero}).
\begin{thm}\label{thm:volume growth rigidity}
We assume $\Ric^{N}_{V} \geq c^{-1}c^2_p\,\kappa_{V,p}\,e^{   -\frac{4(1-\eps)f_{V,p}}{n-1}      }g$,
and $c=(n-1)^{-1}$.
If
\begin{equation}\label{eq:volume growth rigidity assumption}
\varliminf_{r\to +\infty}\frac{\nu_{V,p}(B_{V,r}(p))}{\mathcal{S}_{\kappa}(r)}\geq \omega_{n-1},
\end{equation}
and if $C_{\kappa}=+\infty$,
then $M$ is diffeomorphic to $\mathbb{R}^{n}$,
and the following properties hold:
\begin{enumerate}\setlength{\itemsep}{+1.0mm}
\item If $N=n$, then $V\equiv 0$, and $g=dt^2+\s_{c^2_p\kappa}(t)\,g_{\mathbb{S}^{n-1}}$;
\item if $N=1$, then $\eps=0$, and
\begin{equation*}
g=dt^2+c^{-2}_p\,\exp\left( \frac{2f_{V,p}(\gamma_{v}(t))}{n-1} \right)\,\s^2_{\kappa}(s_{V,v}(t))g_{\mathbb{S}^{n-1}};
\end{equation*}
\item if $N\neq 1,n$, then $\eps=0$, $V$ is orthogonal to $\nabla d_p$ on $M\setminus \{p\}$ and vanishes at $\{p\}$, and $g=dt^2+\s_{c^2_p\kappa}(t)\,g_{\mathbb{S}^{n-1}}$.
\end{enumerate}
\end{thm}
\begin{pf}
Due to Propositions \ref{thm:absolute volume comparison} and \ref{thm:relative volume comparison},
the assumption (\ref{eq:volume growth rigidity assumption}) tells us that
the equality in \eqref{eq:absvolcompa} holds for all $r>0$.
From $C_{\kappa}=+\infty$
we derive $\tau_{V}(v)=+\infty$ for all $v\in U_{p}M$ (see Remark \ref{rem:equality in absolute volume comparison});
in particular,
$\tau(v)=+\infty$,
and $M$ is diffeomorphic to $\mathbb{R}^{n}$.
Now,
the equality in (\ref{eq:Laplacian comparison}) holds on $M\setminus \{p\}$ (see Remark \ref{rem:equality in volume element comparison}).
In virtue of Lemma \ref{lem:LaplacianRigidity},
we complete the proof of Theorem \ref{thm:volume growth rigidity}.
\end{pf}

\begin{remark}\label{rem:refLaplacian}
{\rm The authors do not know whether a similar result holds when $C_{\kappa}<+\infty$.
In this case,
under the same setting as in Theorem \ref{thm:volume growth rigidity},
we see $\tau_{V}(v)=C_{\kappa}$ for all $v\in U_pM $.
Since $\tau(v)$ can be either finite or infinite,
it seems to be difficult to conclude any rigidity results.}
\end{remark}

%%%%%%%%%%%%%%%%%%%%%%%%%%%%
\subsection{Radial case}

Here we consider the case where $f_{V,p}$ is radial with respect to $p$.
In this case,
$s_{V,v}(t)$ does not depend on $v$,
and we can write it as $s_{V}(t)$.
In particular,
Lemma \ref{lem:volume element comparison} can be rewritten as follows:
\begin{lem}\label{lem:radial volume element comparison}
Let $f_{V,p}$ be radial with respect to $p$.
We assume
\begin{equation*}
\Ric^{N}_{V}(\dot{\gamma}_{v}(t))\geq c^{-1}c^2_p\,\kappa(s_{V}(t))\,e^{   -\frac{4(1-\eps)f_{V,p}(\gamma_v(t))}{n-1}      }
\end{equation*}
for all $t \in ]0,\tau(v)[$.
Then for all $t_{1},t_{2} \in ]0,\tau(v)[$ with $t_{1}\leq t_{2}$
\begin{equation*}\label{eq:relative volume element comparison}
\frac{\theta_{V}(t_{2},v)}{ \theta_{V}(t_{1},v)}\leq \frac{\mathfrak{s}^{1/c}_{\kappa}(s_{V}(t_2))}{\mathfrak{s}^{1/c}_{\kappa}(s_{V}(t_1))}.
\end{equation*}
Moreover,
if $c=(n-1)^{-1}$, then
for all $s\in [0,\tau(v)[$ we have
\begin{equation*}\label{eq:absolute volume element comparison}
\theta_{V}(t,v)\leq \mathfrak{s}^{n-1}_{\kappa}(s_{V}(t)).
\end{equation*}
\end{lem}

For $r>0$,
we set
\begin{equation*}
B_{r}(p):=\left\{\,x\in M \mid d_{p}(x) < r \,\right\},\quad \mathcal{S}_{\kappa,V}(r):=\int^{r}_{0}\,\bar{\s}^{1/c}_{\kappa}(s_V(t))\,\d t.
\end{equation*}
Having Lemma \ref{lem:radial volume element comparison} at hand,
we can prove the following assertions along the lines of the proof of the statements in the previous subsections.
The proof is left to the readers.
\begin{prop}
Let $f_{V,p}$ be radial with respect to $p$.
If $\Ric^{N}_{V} \geq c^{-1}c^2_p\,\kappa_{V,p}\,e^{   -\frac{4(1-\eps)f_{V,p}}{n-1}      }g$ and $c=(n-1)^{-1}$,
then for all $r>0$
we have
\begin{equation*}
\mu_{V,p}(B_{r}(p))\leq \omega_{n-1}\,\mathcal{S}_{\kappa,V}(r).
\end{equation*}
In particular,
\begin{equation*}
\varlimsup_{r\to +\infty}\frac{\mu_{V,p}(B_{r}(p))}{\mathcal{S}_{\kappa,V}(r)}\leq \omega_{n-1}.
\end{equation*}
\end{prop}

\begin{prop}
Let $f_{V,p}$ be radial with respect to $p$.
If $\Ric^{N}_{V} \geq c^{-1}c^2_p\,\kappa_{V,p}\,e^{   -\frac{4(1-\eps)f_{V,p}}{n-1}      }g$,
then for all $r,R>0$ with $r\leq R$
we have
\begin{equation*}
\frac{\mu_{V,p}(B_{R}(p))}{\mu_{V,p}(B_{r}(p))} \leq \frac{\mathcal{S}_{\kappa,V}(R)}{\mathcal{S}_{\kappa,V}(r)}.
\end{equation*}
\end{prop}

\begin{thm}
Let $f_{V,p}$ be radial with respect to $p$.
Assume $\Ric^{N}_{V} \geq c^{-1}c^2_p\,\kappa_{V,p}\,e^{   -\frac{4(1-\eps)f_{V,p}}{n-1}      }g$,
and $c=(n-1)^{-1}$.
If we have
\begin{equation*}
\varliminf_{r\to +\infty}\frac{\mu_{V,p}(B_{r}(p))}{\mathcal{S}_{\kappa,V}(r)}\geq \omega_{n-1},
\end{equation*}
and if $C_{\kappa}=+\infty$,
then $M$ is diffeomorphic to $\mathbb{R}^{n}$,
and the following properties hold:
\begin{enumerate}\setlength{\itemsep}{+1.0mm}
\item If $N=n$, then $V\equiv 0$, and $g=dt^2+\s_{c^2_p\kappa}(t)\,g_{\mathbb{S}^{n-1}}$;
\item if $N=1$, then $\eps=0$, and
\begin{equation*}
g=dt^2+c^{-2}_p\,\exp\left( \frac{2f_{V,p}(\gamma_{v}(t))}{n-1} \right)\,\s^2_{\kappa}( s_{V}(t))g_{\mathbb{S}^{n-1}};
\end{equation*}
\item if $N\neq 1,n$, then $\eps=0$, $V$ is orthogonal to $\nabla d_p$ on $M\setminus \{p\}$ and vanishes at $\{p\}$, and $g=dt^2+\s_{c^2_p\kappa}(t)\,g_{\mathbb{S}^{n-1}}$.
\end{enumerate}
\end{thm}

%%%%%%%%%%%%%%%%%%%%%%%%%%%%
%%%%%%%%%%%%%%%%%%%%%%%%%%%%
%%%%%%%%%%%%%%%%%%%%%%%%%%%%
\section{Compactness}\label{sec:Compactness}

%%%%%%%%%%%%%%%%%%%%%%%%%%%%
\subsection{$\eps$-completeness}

We stated that $M$ is compact under the setting of Proposition \ref{thm:finite diameter comparison}.
We first discuss the compactness under that of Proposition \ref{thm:diameter comparison}.
We say that $(M,g,V)$ is \emph{$\varepsilon$-complete at $p$} if 
\begin{align*}
\varlimsup_{r\to+\infty}\inf_{\gamma}\int_0^re^{-\frac{2(1-\varepsilon)f_{V,p}(\gamma(t))}{n-1}}\d t=+\infty,%\label{eq:phimcomplete}
\end{align*}
where the infimum is taken over all unit speed minimal geodesics $\gamma:[0,r]\to M$ with $\gamma(0)=p$ (cf. \cite[Proposition 3.4]{WyYero}, \cite[Definition 2.1]{KL}, \cite[Definition 2.1]{KS}, \cite[Definition 3.2]{LMO:CompaFinsler}).
Note that
in the gradient case of $V=\nabla f$,
the assumption $(1-\eps)f\leq (n-1)\delta$ for $\delta \in \mathbb{R}$ in Proposition \ref{thm:finite diameter comparison} implies the $\eps$-completeness.
We also see the following:
\begin{lem}\label{lem:phimcomplete}
Suppose that 
$(M,g,V)$ is $\varepsilon$-complete at $p$. Then, for any sequence $\{q_i\}$ in $M$  such that  $d_p(q_i)\to+\infty$ as $i\to+\infty$, 
we have $s_p(q_{i})\to+\infty$. 
\end{lem}
\begin{pf}
The proof is similar to that of \cite[Proposition~3.4]{WyYero}. We omit it.   
\end{pf}

Proposition \ref{thm:diameter comparison} together with Lemma \ref{lem:phimcomplete} tells us the following:
\begin{prop}\label{prop:compact}
Let $C_{\kappa}<+\infty$, and assume that $\Ric^{N}_{V} \geq c^{-1}c^2_p\,\kappa_{V,p}\,e^{   -\frac{4(1-\eps)f_{V,p}}{n-1}      }g$, and $(M,g,V)$ is $\varepsilon$-complete at $p$.   Then $M$ is compact. 
\end{prop}

%%%%%%%%%%%%%%%%%%%%%%%%%%%%
\subsection{Ambrose type theorem}

One can further generalize Proposition \ref{prop:compact} in the case where $\kappa$ is a constant function.
Let us prove the following Ambrose type theorem (cf. \cite{Ambrose}):
\begin{thm}\label{thm:AmbroseMyers} 
Assume that $(M,g,V)$ is $\varepsilon$-complete at $p$. 
Suppose additionally that for every unit speed geodesic $\gamma:[0,+ \infty[\to M$ with $\gamma(0)=p$,
we have
\begin{align}
\int_1^{+\infty}e^{\frac{2(1-\varepsilon) f_{V,p}(\gamma(t))}{n-1}}{\Ric}_{V}^N(\dot{\gamma}(t))\d t=+\infty.\label{eq:Ambrose}
\end{align}
Then $M$ is compact. 
\end{thm}
\begin{pf}
Suppose that $M$ is non-compact.
Then there exists a unit speed minimal geodesic $\gamma:[0,+ \infty[\to M$ with $\gamma(0)=p$.
We set
\begin{equation*}
f(t):=f_{V,p}(\gamma(t)),\quad \lambda(t):=e^{\frac{2(1-\eps)f(t)}{n-1}}\,\Delta_{V}d_{p}(\gamma(t)).
\end{equation*}
Note that $\lambda(t)$ is smooth along $\gamma$. 
Lemma \ref{lem:Riccati} leads us to
\begin{align*}
\lambda(r)-\lambda(1)+c\,\int_1^r e^{-\frac{2(1-\varepsilon)f(t)}{n-1}}\lambda(t)^2\d t
&\leq - \int_1^r 
e^{\frac{2(1-\varepsilon)f(t)}{n-1}}\Ric_V^N(\dot{\gamma}(t))\d t.
\end{align*} 
From the assumption \eqref{eq:Ambrose},
\begin{align}
\lim_{r\to+\infty}\left(\lambda(r)+c\,\int_1^r e^{-\frac{2(1-\varepsilon)f(t)}{n-1}}\lambda(t)^2\d t\right)=-\infty.\label{eq:MyerEstimate}
\end{align}
In particular, $\lim_{r\to+\infty}\lambda(r)=-\infty$. 

Next we prove that there exists a finite number $R> 0$ such
that $\lim_{r\uparrow R}\lambda(r)=-\infty$, which contradicts the smoothness of $\lambda(r)$. By \eqref{eq:MyerEstimate}, 
given $C>c^{-1}$ there exists $r_0>1$ such that 
$$
-\lambda(r_0)-c\,\int_1^{r_0}e^{-\frac{2(1-\varepsilon)f(t)}{n-1}}\lambda(t)^2\d t\geq c\,C>1.
$$
By \eqref{eq:Ambrose},
there exists $r_1\in]r_0,+\infty[$ such that 
$\int_{r_0}^r
e^{\frac{2(1-\varepsilon)f(t)}{n-1}} \Ric_V^N(\dot{\gamma}(t))\d t\geq0$ for all $r \geq r_1$. 
Let $\psi(r)$ be the function defined by 
\begin{align*}
\psi(r):=-\lambda(r)-c\,\int_1^r e^{-\frac{2(1-\varepsilon)f(t)}{n-1}}\lambda(t)^2\d t-\int_1^r e^{\frac{2(1-\varepsilon)f(t)}{n-1}} \Ric_V^N(\dot{\gamma}(t))\d t.%\label{eq:MyerEstimateInt}
\end{align*}
By Lemma \ref{lem:Riccati}, we see $\psi'(r)\geq0$. 
Hence $\psi(r)\geq\psi(r_0)$ for $r\geq r_1>r_0$. This implies that 
\begin{align}
-\lambda(r)-c\,\int_1^r e^{-\frac{2(1-\varepsilon)f(t)}{n-1}}\lambda(t)^2\d t\geq c\,C>1\label{eq:Stpe1}
\end{align}
holds for all $r\geq r_1$. 
Let us consider the sequence $\{r_{\ell}\}$ defined inductively by 
\begin{align*}
\int^{r_{\ell+1}}_{r_{\ell}}e^{-\frac{2(1-\varepsilon)f(t)}{n-1}}
\d t
=c^{-1}\left(c\,C\right)^{-\ell+1}\quad \text{ for }\quad\ell\geq1.
\end{align*}
Let $R$ be the increasing limit of $\{r_{\ell}\}$. Then
\begin{align*}
\int_{r_1}^{R}e^{-\frac{2(1-\varepsilon)f(t)}{n-1}}	
\d t=\frac{C}{c\,C-1}<+\infty. 
\end{align*}
In view of the $\varepsilon$-completeness of $(M,g,V)$ at $p$, we obtain $R<+\infty$. 
Finally we claim that for given $\ell \geq 1$, $-\lambda(r)\geq\left( c\,C\right)^{\ell}$ for all $r\geq r_{\ell}$. 
This is true for $\ell=1$ by \eqref{eq:Stpe1}. Suppose that $-\lambda(r)\geq\left( c\,C\right)^{\ell}$ for all $r\geq r_{\ell}$ and fix $r\geq r_{\ell+1}$. Then using inequality \eqref{eq:Stpe1} again, 
\begin{align*}
-\lambda(r)&\geq c\,C+c\,\int_1^{r_{\ell}}
e^{-\frac{2(1-\varepsilon)f(t)}{n-1}}
\lambda(t)^2\d t+ 
c\,\int_{r_{\ell}}^{r}
e^{-\frac{2(1-\varepsilon)f(t)}{n-1}}\lambda(t)^2\d t \\
&\geq  
c\,\int_{r_{\ell}}^{r_{\ell+1}}
e^{-\frac{2(1-\varepsilon)f(t)}{n-1}}
\lambda(t)^2\d t \geq 
c\, (c\,C)^{2\ell}\,\frac{c^{-\ell}}{C^{\ell-1}}=\left(c\,C \right)^{\ell+1}.
\end{align*}
Hence we prove the claim. 
In particular, $\lim_{r\uparrow R}\lambda(r)=-\infty$, which is the desired contradiction.
We complete the proof.
\end{pf}

It is trivial that the condition $\Ric^{N}_{V} \geq c^{-1}c^2_p\,\kappa_{V,p}\,e^{   -\frac{4(1-\eps)f_{V,p}}{n-1}      }g$ implies \eqref{eq:Ambrose} in the case where $\kappa$ is a positive constant and $\eps$-completeness at $p$ holds.

\begin{remark}
{\rm  
In the gradient case of $V=\nabla f$,
the first named author and Li \cite{KL} have proved Theorem \ref{thm:AmbroseMyers} under the curvature condition \eqref{eq:RicciLowerBddKL} (see \cite[Theorem 2.12]{KL}). 
Furthermore,
the first named author and Shukuri \cite{KS} have extended it to the non-gradient case (see \cite[Theorem 2.17]{KS}).}
\end{remark}

%\noindent
%\emph{Acknowledgment.} The authors would 
%thank to 
%the anonymous referee. 
%His/Her comments help to improve the quality 
%of this paper very much.

\providecommand{\bysame}{\leavevmode\hbox to3em{\hrulefill}\thinspace}
\providecommand{\MR}{\relax\ifhmode\unskip\space\fi MR }
\providecommand{\MRhref}[2]{%
  \href{http://www.ams.org/mathscinet-getitem?mr=#1}{#2}
}
\providecommand{\href}[2]{#2}

% \bibliographystyle{amsplain}
% \bibliography{refs}
\end{document}